\documentclass[11pt,a4paper]{article}

\usepackage[latin1]{inputenc}

\usepackage{amsfonts,a4wide}
\usepackage{amsmath,mathenv}

\begin{document}
\numberwithin{equation}{section}
\newtheorem{theoreme}{Theorem}[section]
\newtheorem{proposition}{Proposition}[section]
\newtheorem{remarque}{Remark}[section]
\newtheorem{lemme}{Lemma}[section]
\newtheorem{corollaire}{Corollary}[section]
\newtheorem{definition}{Definition}[section]
\newtheorem{conjecture}{Conjecture}[section]
\renewcommand{\theenumi}{\roman{enumi}}
\def\RR{{\rm I\hspace{-0.50ex}R} }
\def\DD{\displaystyle}
\def\be{\begin{equation} \displaystyle}
\def\ee{\end{equation} }
\def\ben{\[\displaystyle}
\def\een{\] }
\def\bea{\begin{eqnarray}}
\def\eea{\end{eqnarray} }
\def\bean{\begin{eqnarray*}}
\def\eean{\end{eqnarray*}}
\def\div{{\rm div \;}}
\def\NN{{\rm I\hspace{-0.50ex}N} }
\def\ZZ{{\sf Z\hspace{-0.80ex}Z} }
\def\RR{{\rm I\hspace{-0.50ex}R} }
\def\CC{\rm \hbox{C\kern-.57em\raise.47ex
         \hbox{$\scriptscriptstyle |$}\kern+0.5 em }}
\def\erfc{{\rm erfc}}
\newcommand{\dps}{\displaystyle }

\title{Mathematical analysis of a nonlinear parabolic equation arising
  in  the modelling of non-newtonian flows} 
\author{Eric Canc\`es$^{(a)}$, Isabelle Catto$^{(b)}$ and Yousra
  Gati$^{(a)}$ \\
       \footnotesize{(a) CERMICS, Ecole Nationale des
          Ponts et Chauss\'ees and INRIA,} \\ 
           \footnotesize{6 \& 8 avenue Blaise Pascal, Cit\'e
          Descartes, 77455 Marne-la-Vall\'ee Cedex 2, France.} \\
\footnotesize{(b) CEREMADE, UMR CNRS 7534, Universit\'e Paris
 IX-Dauphine,}  \\ \footnotesize{Place du Mar\'echal de Lattre de
 Tassigny, F-75775 Paris 
 Cedex 16, France.}}

\date{\today}
\maketitle

\begin{abstract}
The mathematical properties of a nonlinear parabolic equation arising
  in the modelling of non-newtonian flows are investigated. The
  peculiarity of this equation is that it may degenerate into a
  hyperbolic equation (in fact a linear advection equation). 
Depending on
  the initial data, at least two situations can be encountered: the 
equation
  may have a unique solution in a convenient class, or it may have
  infinitely many solutions. 
\end{abstract}

\section{Introduction}
\label{sec:intro} Modelling the flow of complex fluids is a very
intricate problem which is far from being solved up to now.
Besides studies which aim at improving phenomenological rheological
models (purely macroscopic constitutive laws), only a few attempts are
made to recover the rheological behavior of a complex fluid
from elementary physical processes arising in its microstructure.\\

\medskip

\noindent
The mesoscopic model which has been proposed by H\'ebraud and Lequeux 
in
\cite{bi:HL} deals with simple shear flows of
concentrated suspensions. It is obtained by dividing
the material in a large number of mesoscopic elements 
(``blocks'') with a given shear stress $\sigma$ ($\sigma$ is a real
number; it is in fact an extra-diagonal term of the stress tensor in
convenient coordinates) and by
considering the evolution of the 
probability density $p(t,\sigma)$ which represents the distribution of
stress in the 
assembly of blocks.  Under various assumptions on the evolution of the
stresses of the 
blocks which will be described below, the equation for the 
probability density
$p(t,\sigma)$ for a block to be under stress $\sigma$ at time
$t$   may be written as:
\begin{subequations}
\label{eq:syst-p}
 \begin{EqSystem}
    \partial_{t}p=-b(t) \; \partial_{\sigma}p+D(p(t))\;
    \partial^2_{\sigma\sigma}p-\frac{\chi_{\RR \setminus
    [-\sigma_c,\sigma_c]}
    (\sigma)}{T_0}\;p+\frac{D(p(t))}{\alpha}\;\delta_{0}(\sigma)\quad
    \textrm{ on } 
(0;T)\times \RR\;;\label{eq:p}\\
p\,\ge\,0\;;\\
p(0,\sigma)=p_0(\sigma)\;, \label{eq:p0er}
  \end{EqSystem}
\end{subequations}
 where for $f \in L^1(\RR)$, we denote by
\begin{equation*}
D(f) = \frac{\alpha}{T_0} \int_{|\sigma|>\sigma_{c}} f(\sigma)\, 
d\sigma\;.
\end{equation*}
In equation~(\ref{eq:p}), $\chi_{\RR \setminus [-\sigma_c,\sigma_c]}$
  denotes the characteristic function of the open set $\RR \setminus
  [-\sigma_c,\sigma_c]$ and $\delta_0$ the Dirac delta function on
  $\RR$. Each term arising in the above equation \eqref{eq:p} (HL 
equation
in short) has a clear physical interpretation. When a block is
sheared, the stress of this block evolves with a
variation rate $b(t) = G_0 \dot \gamma(t)$ proportional to the shear
rate $\dot \gamma(t)$ ($G_0$ is an
elasticity constant); in this study, the shear rate $\dot \gamma(t)$, 
and
therefore the function $b(t)$, are assumed to be in $L^2_{\rm
  loc}(\RR^+)$. When the modulus of the stress  
overcomes a critical value $\sigma_c$, the block becomes unstable
and may relax into a state with zero stress after a characteristic 
relaxation
time $T_0$. This phenomenon induces a rearrangement
of the blocks and is modelled through the diffusion term
$D(p(t))\;\partial^2_{\sigma\sigma}p$. The diffusion coefficient 
$D(p(t))$ is assumed to be proportional to the amount of stress
which has to be redistributed by time unit and the positive
parameter $\alpha$ is supposed to represent the ``mechanical
fragility" of the material. 

\medskip

\noindent
In all that follows,  the parameters $\alpha$, $T_0$ and $\sigma_c$  
are
positive, and the initial data $p_0$ in \eqref{eq:p0er} 
is  a given probability density; that is
\begin{equation}\label{IC}
 p_0\,\ge\,0\,,\quad p_0\in L^1(\RR)\,, \quad
 \int_{\RR}p_0 =1\;.
\end{equation}
We will be looking for  solutions $p=p\,(t,\sigma)$ in
$C^0_{t}(L^1_\sigma\cap L^2_\sigma)$ such that  $\sigma\,p$ belongs to
$\displaystyle{ L^\infty_t(L^1_\sigma)}$  to the nonlinear
parabolic partial differential equation~\eqref{eq:syst-p}.
The subscript $\sigma$ refers to integration over $\RR$ with respect 
to
$\sigma$, whereas the subscript $t$ refers to time integration on
$[0,T]$ for any $T > 0$. 
Note that the average stress in the material is given by
\begin{equation} \label{eq:def-tau}
\tau(t)=\int_\RR\sigma\,p(t,\sigma)\,d\sigma \;,
\end{equation}
and therefore the above condition on $\sigma p$ ensures that the
average stress is an essentially bounded  function of  time.\\

\medskip

\noindent
Actually in practice, the shear rate is not uniform in the flow and 
in order to
better describe the coupling  
between the macroscopic flow and the evolution of the microstructure 
we
introduce and study in a 
second  paper~\cite{bi:CCGL} a micro-macro model. In this model  $p$ 
is  also a function of the macroscopic 
space variables and  the 
average stress defined by~\eqref{eq:def-tau} is inserted into the 
macroscopic 
equation governing the velocity of the macroscopic flow.

\medskip

\noindent
In order to lighten the notation and without loss of
generality  we  assume from now on that
$\sigma_c=1$ and $T_0=1$. This amounts to changing the time and stress
scales. 

\medskip

\noindent
The main difficulties one encounters in the mathematical analysis  
come
from the nonlinearity  
in the diffusion term and also and  even  more from the fact that the
parabolic equation may degenerate  
when the viscosity coefficient $D(p)$ vanishes, and this will be shown
to may appear only when 
$D(p_0)=0$. This difficulty is illustrated on a  simplified example 
just
below and also in 
Section~\ref{sec:stat} where we discuss the existence of
stationary solutions  in the 
case when the shear rate $b$ is a constant. 

\medskip

\noindent
Let us  first of all look at the following simplified model which
already includes  
the difficulties we are going to face to in the study of equation 
~\eqref{eq:syst-p}. We consider the equation~:
\begin{equation}\label{eq:simpl}
\left\{\begin{array}{rlc}
\partial_t u&=&D(u(t))\,\partial^2_{\sigma\sigma}u\;;\\
u(0, \sigma)&=&\frac 1 2\chi_{]-1,1[}(\sigma)\;,
\end{array}\right.
\end{equation}
where $\chi_{]-1,1[}$ is the characteristic function of the
interval $]-1,1[$. The initial 
condition is on purpose chosen in such a way that
$D(u(t=0))=0$. The function $u=\frac 1 2\chi_{]-1,1[}(\sigma)$ is a
stationary solution to this equation and for this solution $D(u(t))$ 
is
identically zero.  But it is not the unique solution 
to~\eqref{eq:simpl} in
$C^0_t(L^2_\sigma)\cap L^\infty_t(L^1_\sigma)$. It is indeed possible 
to
construct a so-called \textit{vanishing viscosity solution} for which
$D(u(t))>0$ for all $t>0$,   and there are  actually infinitely many
solutions to this 
equation. (This statement is obtained as a corollary of
Lemma~\ref{lem:2} in Section~\ref{sec:deg} below.) 

\medskip

\noindent
As far as equation (\ref{eq:syst-p}) is concerned, we show that, in 
the
case when $D(p_0)=0$ and $b \equiv 0$, we may have either a unique or
infinitely many solutions, depending on the initial data (see
Proposition~\ref{prop:deg} in Section~\ref{sec:deg}).

\medskip

\noindent
On the other hand, we are able to prove the following existence and 
uniqueness
result in the non-degenerate case when $D(p_0) > 0$~:
\begin{theoreme}\label{th:main1}
 Let the initial data $p_0$ satisfy the conditions
\begin{equation}\label{eq:IC-full}
p_0\in L^1(\RR)\cap L^\infty(\RR)\,,\quad p_0\geq 0\,,\quad \int_\RR
p_0=1\quad\mathrm{ and }\int_\RR|\sigma|\,p_0<+\infty\;,
\end{equation}
and assume that 
\begin{equation*}
D(p_0)>0\;.
\end{equation*}
Then, for every $T>0$, there exists a
unique solution $p$ to the system~\eqref{eq:syst-p} in
$L^\infty_t(L^1_\sigma \cap L^2_\sigma)\cap L^2_t(H^1_\sigma)$. 
Moreover, $p\in L^\infty_{t,\sigma}\cap C^0_t(L^1_\sigma\cap 
L^2_\sigma)$,
$\int_\RR p(t,\sigma) \, d\sigma=1$ for all $t > 0$,
$D(p)\in C^0_t$ and for every 
$T>0$ there exists a positive constant $\nu(T)$ such that
\begin{equation*}
\min_{0\leq t\leq T}D(p(t))\geq \nu(T)\;.
\end{equation*}
Besides $\sigma\,p\in L^\infty_t(L^1_\sigma)$ so that the average 
stress
$\tau(t)$ is well-defined by (\ref{eq:def-tau}) in $L^\infty_t$.
\end{theoreme}

\medskip

\noindent
The first  step toward the existence proof of solutions
to~\eqref{eq:syst-p} will consist in 
the study of  so-called vanishing viscosity approximations, which are
the unique solutions to the family of equations
\begin{subequations}\label{syst:p-eps}
\begin{EqSystem}
\partial_tp_\varepsilon\,=\,-\,b(t)\,\partial_\sigma
p_\varepsilon\,+\,(D(p_\varepsilon(t))+\varepsilon)\,\partial^2_{\sigma\sigma}p_\varepsilon\,
-\,\chi_{\RR\setminus [-1,1]}\,p_\varepsilon\,+
\frac{D(p_\varepsilon(t))}{\alpha}\,\delta_0(\sigma)\; ;\label{eq:p-eps}\\
p_\varepsilon\geq 0\;;\\
p_\varepsilon(0,\cdot)\,=p_0\;
  \end{EqSystem}
\end{subequations}
(recall that we have rescaled the time and stress units to get $T_0=1$
and $\sigma_c=1$).
Section \ref{sec:visc} below is devoted to the proof of the following
\vskip10pt
\begin{proposition}[Existence and uniqueness of vanishing viscosity
approximations]\label{prop:visc} Let $T>0$ be given. We assume that 
the initial data satisfies the same conditions
\eqref{eq:IC-full} as in the statement of the theorem. Then, for 
every $T>0$ and $0<\varepsilon\leq 1$, there exists a unique solution
$p_\varepsilon$  to \eqref{syst:p-eps} in $L^\infty_t(L^1_\sigma\cap
L^2_\sigma)\cap L^2_t(H^1_\sigma)$.
Moreover,  $\,p_\varepsilon\in L^\infty_{t,\sigma}\cap
C^0_t(L^1_\sigma\cap L^2_{\sigma})\,$, $D(p_\varepsilon)\in C^0_t$,
\begin{equation}\label{eq:proba-eps}
\int_\RR p_\varepsilon=1\;,
\end{equation}
\begin{equation}\label{eq:bd-infini-eps}
0\leq p_\varepsilon\leq \Vert
p_0\Vert_{L^\infty_\sigma}+\sqrt{\frac{\alpha}{\pi}}\,\sqrt{T}\;,
\end{equation}
and for every $T>0$, there exist positive constants $C_1(T,p_0)$, 
$C_2(T,p_0)$ and $C_3(T,p_0)$ which
are independent of
$\varepsilon$  such that
\begin{equation}\label{eq:bd-sigma-eps}
\sup_{0\leq t\leq T}\int_\RR |\sigma| \, p_\varepsilon\leq C_1(T,p_0)
\,,
\end{equation}
\begin{equation}\label{eq:bd-L2-eps}
\sup_{0\leq t\leq T}\int_\RR p_\varepsilon^2\leq C_2(T,p_0)\;,
\end{equation}
and
\begin{equation}\label{eq:bd-grad-eps}
\int_0^T\big(\varepsilon+D(p_\varepsilon)\big)\,\int_{\RR}\vert
  \partial_{\sigma}p_\varepsilon \vert^2
\,\le\,C_3(T,p_0)\;.
\end{equation}
\end{proposition}

\medskip

\noindent
Theorem (\ref{th:main1}) is then proved in Section~\ref{sec:nondeg}
while the degenerate case is investigated in Section~\ref{sec:deg}. 
Lastly, the description of stationary solutions in the constant shear
rate case is carried out in Section~\ref{sec:stat}. 

\vskip10pt
\section{The vanishing viscosity approximation}\label{sec:visc}
This section is devoted to the proof of Proposition~\ref{prop:visc}. 
\\

\medskip

\noindent
We begin with the following~:
\begin{lemme}[Uniqueness]\label{lem:unique}
Let $p_0$ satisfy \eqref{IC}. Then for every $T>0$ and 
$0<\varepsilon$, there exists at most one
solution $p_\varepsilon$ to \eqref{syst:p-eps} in 
$L^\infty_t(L^1_\sigma\cap L^2_\sigma)\cap
L^2_t(H^1_\sigma)$. Moreover, $p_\varepsilon\in C^0_t(L^2_\sigma)$ 
(thus,
the initial condition makes sense) and
\begin{equation}\label{eq:proba}
\int_\RR p_\varepsilon =1\;,
\end{equation}
for almost every $t$ in $[0,T]$.
\end{lemme}
\vskip10pt\noindent\textbf{Proof of Lemma~\ref{lem:unique}:} We begin 
with proving that every solution to
\eqref{syst:p-eps} in $L^\infty_t(L^1_\sigma\cap L^2_\sigma)\cap
L^2_t(H^1_\sigma)$ satisfies
\eqref{eq:proba}. We fix $R>1$ and we consider a cut-off $C^2$ 
function
$\phi_R=\phi_R(\sigma)$ with compact support which is  equal to 1 when
$\displaystyle{0\,\le\,\vert\sigma\vert\,\le\,R}$ and to 0 when
$\displaystyle{\vert\sigma\vert\,\ge\,2R}$ and such that
\begin{equation}
\label{C.1}
\vert\phi'_R\vert\le\frac{C}{R}\;,
\end{equation}
where here and below $C$ denotes a  positive constant that is 
independent of $R$. Notice
that
$\phi '$ is equal to 0 on $]-\infty,\,-2R]$, on $[-R,\,R]$ and on
$[2R,\,+\infty[$.

\medskip

\noindent
Now, we multiply  (\ref{eq:p-eps}) by $\phi _R$ and integrate over
$[0,t]\times \RR$ to obtain
\begin{eqnarray*}
\int_{\RR}p_\varepsilon(t)\,\phi _R-\,\int_{\RR}p_0\,\phi _R&=&
-\int_0^t\,b(s)\,\int_\RR\partial_\sigma p_\varepsilon(s)
\,\phi_R\,
-\,\int_0^t\big(D(p_\varepsilon(s))+\varepsilon\big)\,\int_\RR\partial_{\sigma}p_\varepsilon(s)\,\phi_R'\\ & &-\, 
\int_0^t\int_{\vert\sigma\vert>1}p_\varepsilon(s)\,\phi_R\,
+\,\frac{1}{\alpha}\int_0^tD(p_\varepsilon(s))\,\phi_R(0)\;.
\end{eqnarray*}
We bound from above  the terms on the right-hand side as follows.
First, we have
\begin{eqnarray*}
\Big\vert\int_0^t b(s) \int_\RR \partial_\sigma p_\varepsilon(s)\,
\phi_R\Big\vert&\le&\int_0^t
\,\vert\,b(s)\,\vert\int_\RR\,p_\varepsilon(s)\;\vert
\,\phi'_R\,\vert\\
&\le&\frac{C}{R}\int_0^t
\,\vert\,b(s)\,\vert\,\int_{R\leq \vert\sigma\vert\leq 2R}
p_\varepsilon(s)\,\le\,\frac{C}{R}\;,
\end{eqnarray*}
thanks to (\ref{C.1}) and using that $p_\varepsilon\in
 L^\infty_{t}(L^1_\sigma)$ and
$b\in L^1_{t}$. Next,
\begin{eqnarray*}
\int_0^t\big(D(p_\varepsilon)+\varepsilon\big)\,
\vert\int_\RR\partial_{\sigma}p_\varepsilon\,\phi'\vert
 &\le& (\varepsilon+\alpha\,\Vert
p_\varepsilon\Vert_{L^\infty_t(L^1_\sigma)})\,\int_0^t\Vert\partial_\sigma
p_\varepsilon\Vert_{L^{2}_{\sigma}}\,\Vert\phi'_R\Vert_{L^2_\sigma}\\
&\le&\frac{C\sqrt{t}}{R^{1/2}}\,
\,\Vert\partial_\sigma
p_\varepsilon\Vert_{L^2_{t,\sigma}}\leq \frac{C}{R^{1/2}}\;,
\end{eqnarray*}
thanks again to (\ref{C.1}), Cauchy-Schwarz' inequality and
since   $\partial_\sigma p_\varepsilon$ is in $L^2_{t,\sigma}$. 
Finally,
\begin{eqnarray*}
0\leq\frac{1}{\alpha}\int_0^tD(p_\varepsilon)\,-\,
\int_0^t\int_{\vert\sigma\vert>1}p_\varepsilon\,\phi_R &=&
\int_0^t\int_{\vert\sigma\vert>1}p_\varepsilon\,(\,1\,-\,\phi_R\,)\\ 
&\le&
\int_0^t\int_{\vert\sigma\vert>R}p_\varepsilon\;,
\end{eqnarray*}
and the right-hand side goes to $0$ as $R$ goes to infinity since $ 
p_\varepsilon$ is in
$L^\infty_{t}(L^{1}_{\sigma})$.
 All this together yields
\begin{equation*}
\int_\RR p_\varepsilon(t)=\lim_{R\rightarrow
 +\infty}\int_{\RR}p_\varepsilon(t)\,\phi_R\,=\lim_{R\rightarrow
 +\infty}\int_{\RR}p_0\,\phi_R=\int_{\RR}p_0=1\;,
\end{equation*}
for almost every $t$ in $[0,T]$. In particular, this implies that $ 
D(p_\varepsilon)\leq
\;\alpha\;$.\\

\medskip

\noindent
Let us now argue by contradiction by assuming that there exist two
solutions $p_1$ and $p_2$ to (\ref{eq:p-eps}) corresponding to the 
same
initial data $p_0$. By subtracting  the
equations satisfied by $p_1$ and $p_2$ respectively, we  obtain
  \begin{EqSystem}
\DD{\partial_t q\,=\,-\,b(t)\;\partial_\sigma
q\,+\,D(q)
\;\partial^2_{\sigma\sigma}p_1\,+\,(D(p_2)+\varepsilon)\,\partial^2_{\sigma\sigma}q
-\,\chi_{\RR\setminus 
[-1,1]}\,q\,+\,\frac{D(q)}{\alpha}\,\delta_0(\sigma)}\;;\label{eq:q}\\
\DD{q(0,\sigma)\,=\,0\;,}\nonumber
  \end{EqSystem}
where $q=p_1-p_2$. We multiply (\ref{eq:q}) by $q$ and integrate over 
$\RR$
with respect to $\sigma$ to obtain, after integrations by parts,
\begin{equation}\label{eq:step}
\frac 1 2 \frac{d}{dt} \int_\RR
q^2\,+\,(D(p_2)+\varepsilon)\,\int_\RR\vert\partial_{\sigma}q\vert^2\,+\int_{\vert\sigma\vert>1}q^2\,
=\,\frac{D(q)} \alpha
\,q(t,0)-D(q)\,\int_\RR\partial_\sigma  p_1\partial_\sigma q\;.
\end{equation}
We first remark that since $\DD{\int_\RR p_1 \,=\,\int_\RR 
p_2\,=\,1}$ thanks to
\eqref{eq:proba}, we get
\begin{equation*}
\vert D(q)\vert
=\alpha\,\big\vert \int_{\vert \sigma\vert
  <1}q\big\vert\leq \alpha\sqrt{2}\,\Vert q\Vert_{L^2_\sigma}\;,
\end{equation*}
with the help of Cauchy-Schwarz' inequality. Next, using the Sobolev 
embedding of $H^1(\RR)$ into $L^\infty(\RR)$, we bound from above the 
terms on the right-hand side in the following way~:
\begin{eqnarray*}
\lefteqn{\vert \,\frac{D(q)} \alpha
\,q(t,0)-D(q)\,\int_\RR\partial_\sigma  p_1\partial_\sigma
q\,\vert}\\
&\le&\sqrt{2}\,\Vert q\Vert_{L^2_\sigma}\,\Vert 
q\Vert_{L^\infty_\sigma}\,+
\,\sqrt{2}\,\alpha\,\Vert q\Vert_{L^2_\sigma}\int_\RR\vert 
\partial_\sigma p_1 \partial_\sigma q\vert\\
&\le&\sqrt{2}\,\Vert q\Vert_{L^2_\sigma}\, \big(\,\Vert
q\Vert^2_{L^2_{\sigma}}\,+\,\Vert\partial_\sigma
q\Vert^2_{L^2_\sigma}\,\big)^{\frac 1 2}\,+\,\sqrt{2}\,\alpha\,\Vert 
q\Vert_{L^2_\sigma}\, \Vert
\partial_\sigma p_1\Vert_{L^2_\sigma}\Vert\partial_\sigma 
q\Vert_{L^2_\sigma}\\
&\le& \frac{1}{\varepsilon}\,\Vert
q\Vert^2_{L^2_\sigma}\,+\,\frac {\alpha^2}
{\varepsilon} \Vert
q\Vert^2_{L^2_\sigma}\,\Vert
\partial_\sigma p_1\Vert^2_{L^2_\sigma}\,+\,\frac \varepsilon 2 \Vert
q\Vert^2_{L^2_\sigma}\,+\,\varepsilon \Vert \partial_\sigma
q\Vert^2_{L^2_\sigma}\;.
\end{eqnarray*}
 Therefore,
comparing with \eqref{eq:step} we deduce
\begin{eqnarray*}
\frac 1 2 \frac{d}{dt} \Vert q\Vert^2_{L^2_\sigma}\,\le\,\big(\frac
{1}{\varepsilon}\,+\,\frac{\alpha^2}{\varepsilon}\,\Vert
\partial_\sigma p_1\Vert^{2}_{L^2_\sigma}\,+\,\frac \varepsilon 
2\,\big)\Vert q\Vert^2_{L^2_\sigma}\;.
\end{eqnarray*}
Finally, by applying the Gronwall lemma, we prove that $\DD{\Vert
  q\Vert^2_{L^2_\sigma}\,\le\,0}$, thus $q=0$. The uniqueness of the
  solution follows.
\hfill $\diamondsuit$\\

\begin{remarque} \label{rem:unic}The same proof shows that if there 
exists a solution to \eqref{eq:syst-p} in 
$L^\infty_t(L^1_\sigma \cap L^2_\sigma) \cap  L^2_t(H^1_\sigma)$ such
that $\inf_{0\leq t\leq T}D(p(t))> 0$, then it is unique in this 
space. 
\end{remarque}
{}We now turn to the existence part in the statement of
Proposition~\ref{prop:visc}. {}From now on we fix a positive
constant  $\varepsilon\leq 1$. The proof of
Proposition~\ref{prop:visc} will be carried out by the Schauder
fixed point theorem. For given positive constants $M(\geq
\varepsilon)$ and $R$, we introduce ${\cal D}_{\varepsilon,M}$ and
$Y_R$ two closed convex subsets of respectively $\DD{L^2_t}$ and
$\DD{L^2_{t,\sigma}}$ as follows:
\begin{eqnarray*}
{\cal D}_{\varepsilon,M}\,&=&\,\{a\,\in\,L^2_t;\,\varepsilon\le a\le 
M\}\\
Y_{R}\,&=&\,\{p\in\,L^2_{t,\sigma};\;\,p\ge\,0\;,\; \sup_{0\leq t\leq 
T}\int_\RR \vert \sigma\vert
\,p\,\le\,R\,\}.
\end{eqnarray*}
To simplify notation we denote
\begin{equation*}
\left\{
\begin{array}{rcl}
\varphi_\eta(x)&=&\displaystyle{\frac{1}{\sqrt{2\pi}\;\eta}}\exp\big(-\frac{x^2}{2\,\eta^2}\big)
\quad \mathrm{ if \;}\eta>0\;;\\
\varphi_0&=&\delta_{0} \;.
\end{array}
\right.
\end{equation*}

\medskip

\noindent
We first prove the following
\begin{proposition}\label{prop:1} Let
$T>0$ and let $p_0\in L^2(\RR) $ such that $p_0\geq 0$. Then, for 
every   $a$ in ${\cal
D}_{\varepsilon,M}\,$ and $q$ in $Y_R$,  there exists a unique  
solution  $p$ in
$L^\infty_t(L^2_\sigma)\,\cap\,L^2_t(H^1_\sigma)$ to
 \begin{subequations}\label{eq:syst-lin}
 \begin{EqSystem}
\partial_tp(t,\sigma)\,=\,-\,b(t)\,\partial_\sigma
p(t,\sigma)\,+\,a(t)\,\partial^2_{\sigma\sigma}p(t,\sigma)\,-\,\chi_{\RR\setminus 
[-1,1]}(\sigma)
\,p(t,\sigma)
\,+
\frac{D(q)}{\alpha}\,\delta_0(\sigma)\; ;\label{eq:plin}\\
p(0,\sigma)\,=p_0(\sigma)\;.
  \end{EqSystem}
\end{subequations}
Moreover, $p\in C^0_t(L^2_\sigma)$, $p$ is non-negative  and
\begin{equation}\label{eq:compare}
p_-\leq p\leq p_+\;,
\end{equation}
with
\begin{equation}\label{eq:p-}
p_{-}(t,\sigma)= e^{ -t }
\int_{-\infty}^{+\infty} p_{0}(\sigma')\,
\varphi_{\sqrt{2\int_0^ta\,}}(\sigma-\sigma'-\chi(t) )\,d\sigma'
\end{equation}
and
\begin{eqnarray}\label{eq:p+}
p_{+}(t,\sigma)&=&
\int_{-\infty}^{+\infty} p_{0}(\sigma')\,
\varphi_{\sqrt{2\int_0^ta\,}}(\sigma-\sigma'-\chi(t) 
)\,d\sigma'+\nonumber\\
& &+\frac{1}{\alpha}\int_0^t D(q(s))\,
\varphi_{\sqrt{2\int_s^ta\,}}(\sigma-\chi(t) +\chi(s))\,ds\;,
\end{eqnarray}
where  $\chi(t)=\int_0^t b(s)\,ds$. In addition,
\begin{enumerate}
\item If  $p_0\in L^\infty(\RR)$,  then  $p$ is in 
$L^\infty_{t,\sigma}$ and
\begin{equation}\label{eq:bd-infini}
0\leq p\leq  \Vert
p_0\Vert_{L^\infty}+\frac{R\;\sqrt{T}}{\sqrt\pi\,\sqrt\varepsilon}\;.
\end{equation}
\item If $\int_\RR |\sigma|\,p_0<+\infty\,$ (\/thus
$p_0\in L^1(\RR)$\/),
then $\vert\sigma|\,p\in L^\infty_t(L^1_\sigma)$. More precisely, we 
have
\begin{eqnarray}\label{eq:bd-sigmap}
\sup_{0\leq t\leq T}\int_\RR|\sigma|\,p
&\leq & \int_\RR |\sigma|\,p_0+\sqrt{T}\,\Vert
b\Vert_{L^2(0,T)}\;\Vert p_0\Vert_{L^1}+\frac{2\,R}{3}T^{3/2}\,\Vert
b\Vert_{L^2(0,T)}\nonumber\\
&&+\frac{2}{\sqrt{\pi}}(M\,T\big)^{1/2}\;\Vert
p_0\Vert_{L^1}+\frac{4\,R\,\sqrt{M}}{3\,\sqrt\pi}
\,T^{3/2}\;.
\end{eqnarray}
Moreover $p\in C^0_t(L^1_\sigma)$ and  $D(p)\in
C^0_t$.
\end{enumerate}
\end{proposition}
\vskip6pt\noindent\textbf{Proof of Proposition~\ref{prop:1}:} Let us 
first observe that for  every
$q$ in $Y_R$, $\,D(q)\in L^\infty_t$ since
\begin{equation}\label{eq:bd-D}
0\leq D(q(t))\leq \alpha \int_{\vert \sigma\vert>1}\vert 
\sigma\vert\, q\leq
\alpha\,R\;,
\end{equation}
for almost every $t$ in $[0,T]$. Therefore the source term 
$D(q(t))\,\delta_0(\sigma)$ in
\eqref{eq:plin} is in
$L^\infty_t (H^{-1}_\sigma)$ and the existence and the  uniqueness of 
a solution $p\in
C^0_t(L^2_\sigma)\cap L^2_t(H^1_\sigma)$ to the system
\eqref{eq:syst-lin} is well-known (see for example \cite{bi:chi}). In
particular, 
the initial condition makes sense. Owing to the fact that the source 
term is non -negative, the proof that $p\geq 0$ is also standard  (see again
\cite{bi:chi}).
\medskip

\noindent
We now check the pointwise inequality
\eqref{eq:compare}.\\
This is ensured by the maximum principle with observing that  $p_{-}$ 
and $p_+$ given respectively  by
\eqref{eq:p-} and \eqref{eq:p+} are the unique  solutions to the  
systems
\begin{equation}\label{eq:syst-lin-}
  \left\{
\begin{array}{rcl}
\partial_{t} p_{-}&=&-\,b\,
  \partial_{\sigma}
p_{-}+\,a\,\partial_{\sigma\sigma}^{2}
  p_{-}\,-\,p_{-}\;;\\
p_{-}(0,\sigma)&=&p_0(\sigma)\;,
\end{array}
\right.
\end{equation}
and
\begin{equation}\label{eq:syst-lin+}
  \left\{
\begin{array}{rcl}
\partial_{t} p_{+}&=&\displaystyle{-\,b\,
  \partial_{\sigma}
p_{+}+\,a\,\partial_{\sigma\sigma}^{2}
  p_{+}+\frac{D(q)}{\alpha}\,\delta_0(\sigma)}\;;\\
p_{+}(0,\sigma)&=&p_0(\sigma)\;,
\end{array}
\right.
\end{equation}
respectively. We  now turn to the proof of statement \textit{i.} and
assume that $p_0$  belongs to $L^\infty(\RR)$. Then, using the
two facts that for every
$\nu>0$, 
$\, \int_\RR \varphi_\nu=1$  and $\varphi_\nu\leq 
\frac{1}{\sqrt{2\pi}\;\nu}$, \eqref{eq:bd-infini}
is easily deduced  from $p\leq p_+$ with the help of \eqref{eq:bd-D} 
and
since $a\geq \varepsilon$.\\

\medskip

\noindent
Suppose now that  $\int_\RR |\sigma|\, p_0<+\infty$. This together 
with
the assumption  
$p_0 \in L^2(\RR)$, guarantees that $p_0\in L^1(\RR)$ (see
also below). Using 
\eqref{eq:compare} again, we now have
\begin{eqnarray}
\int_\RR|\sigma|\,p&\leq& \int_\RR|\sigma|\,p_+\nonumber\\
&\leq& \int_\RR \int_\RR
p_0(\sigma')\,|\sigma|\,\varphi_{\sqrt{2\int_0^ta\,}}(\sigma-\chi(t)
-\sigma')\,d\sigma d\sigma' \nonumber\\ & & +\frac{1}{\alpha}\int_0^t
D(q(s))\,\Big(\int_\RR
|\sigma|\,\varphi_{\sqrt{2\int_s^ta\,}}(\sigma-\chi(t)+\chi(s))\,d\sigma
\Big)\,ds\nonumber\\ &=& \int_\RR \int_\RR
p_0(\sigma')\,|\sigma+\sigma'+\chi(t)|\,\varphi_{\sqrt{2\int_0^ta\,}}(\sigma)\,d\sigma
d\sigma' \nonumber\\
 & & +\frac{1}{\alpha}\int_0^tD(q(s))\,\Big(\int_\RR
|\sigma+(\chi(t)-\chi(s))|\,\varphi_{\sqrt{2\int_s^ta\,}}(\sigma)\,d\sigma
\Big)\,ds\nonumber\\
&\leq & \int_\RR |\sigma|\,p_0(\sigma)\,d\sigma+|\chi(t)|\,\Vert 
p_0\Vert_{L^1}+\frac 1\alpha
\int_0^t|\chi(t)-\chi(s)|\,
D(q(s))\,ds\nonumber\\
&&+\frac{2}{\sqrt{\pi}}\big(\int_0^ta\big)^{1/2}\;\Vert 
p_0\Vert_{L^1}+\frac 2{\alpha\sqrt\pi}
\int_0^t 
D(q(s))\, \big(\int_s^ta\big)^{1/2}\,ds\;,\label{eq:bd-sigma}
\end{eqnarray}
since $\int_\RR 
|\sigma|\,\varphi_\nu(\sigma)\,d\sigma=(2/\pi)^{1/2}\,\nu$ and
$\int_\RR\varphi_\nu=1$. With the help of \eqref{eq:bd-D} and
observing that $|\chi(t)-\chi(s)|\leq
\sqrt{t-s}\;\Vert b\Vert_{L^2(0,T)}$, we then deduce 
\eqref{eq:bd-sigmap}. \\

\medskip

\noindent
We now use this bound to check that $p\in C^0_t(L^1_\sigma)$ and
$D(p)\in C^0_t$. Indeed, for any $t$, any sequence $t_n$ in $[0,T]$ 
which converges to $t$  and $A>1$,
we have
\begin{eqnarray}
\int_\RR\big\vert
p(t_n)-p(t)\big\vert&=&\int_{\vert\sigma\vert\leq A}\big\vert
p(t_n)-p(t)\big\vert+\int_{\vert\sigma\vert\geq A}\big\vert
p(t_n)-p(t)\big\vert\nonumber\\
&\leq & \sqrt{2\,A}\,\Big(\int_\RR\big\vert
p(t_n)-p(t)\big\vert^2\Big)^{1/2}+
\frac{1}A\,\int_\RR\vert\sigma\vert\,\big(|p(t_n)|+|p(t)|\big)\nonumber\\
&\leq & \sqrt{2\,A}\,\Big(\int_\RR\big\vert
p(t_n)-p(t)\big\vert^2\Big)^{1/2}+
\frac{2}A\,\sup_{0\leq t\leq T}\int_\RR\vert\sigma\vert\,|p(t)|\;.
\label{eq:cont-L1}\end{eqnarray}
{}For any fixed $A$ the first term in the right-hand side goes to $0$ 
as $n$ goes to infinity  since
$p\in C^0_t(L^2_\sigma)$ and then the second term is arbitrarily 
small as $A$ goes to
infinity. The same argument yields the continuity of $D(p(t))$ with 
respect to $t$.

\hfill$\diamondsuit$\\

\medskip

\noindent
The following proposition aims at checking the required assumptions 
to apply the Schauder fixed point theorem.
\vskip10pt
\begin{proposition}\label{prop:2}  Let $T_f>0$ be given. We assume 
that
\begin{equation}\label{eq:IC}
p_0\in L^1(\RR)\cap L^\infty(\RR)\,,\quad p_0\geq 0\,,\quad \int_\RR
p_0=1\quad\mathrm{ and }\int_\RR|\sigma|\,p_0<+\infty\;.
\end{equation}
Let $0<\varepsilon\leq 1$, $R=1+\int_\RR
|\sigma|\, p_0$ and
$M= 1+ 2\,\alpha$. We define
\begin{equation}\label{eq:def-tc}
T_c=\frac{9}{25}\,\Big[\,\Vert
b\Vert_{L^2(0,T_f)}+\frac{2\,\sqrt{1+2\alpha}}{\sqrt{\pi}}
\Big]^{-\,2}\;.
\end{equation}
Then, for every $T\leq \min\big(\frac{1}{R};T_c\big)$, the
function ${\cal T}:(a;q)\mapsto \big(D(p)+\varepsilon\,; p\big)$,
with $p$ being the solution to the system \eqref{eq:syst-lin},
maps ${\cal D}_{\varepsilon,M}\times Y_R$ into itself. Moreover 
${\cal T}$ is
continuous and ${\cal T}({\cal D}_{\varepsilon,M}\times Y_R)$ is
relatively compact in $L^2(0,T)\times L^2_{t,\sigma}$.
\end{proposition}
\vskip6pt\noindent\textbf{Proof of Proposition~\ref{prop:2}:}
\medskip

\noindent\textbf{Step 1:} \textit{${\cal T}$ is well-defined.}

\medskip

\noindent
According to Proposition~\ref{prop:1},  $p$ is in $C^0_t(L^1_\sigma)$ 
and
$D(p)\in C^0_t$. We now prove that  with our choice for $M$ (
which ensures that
$\varepsilon+D(p_0)\leq 1+\alpha\leq M$\/),
$D(p)+\varepsilon \,\in\,\mathcal{D}_{\varepsilon, M}\;.$ For this, 
we again use the inequality
$p\leq p_+$, the definition
\eqref{eq:p+} of $p_+$, the rough estimate 
$\int_{|\sigma|>1}\varphi_\nu\leq\int_\RR
\varphi_\nu=1$ and
\eqref{eq:bd-D} to obtain
\begin{equation*}
\sup_{0\leq t\leq T}D(p(t))\leq \sup_{0\leq t\leq T} D(p_+(t))\leq
\alpha+\alpha\,R\,T\leq 2\,\alpha\;,
\end{equation*}
for $T\leq \frac{1}{R}$. It only remains now to check that 
$\dps\sup_{0\leq t\leq
T}\int_\RR |\sigma|\,p\leq R$. We thus go back to 
\eqref{eq:bd-sigmap} and observe that this
condition holds provided
\begin{equation*}
T\leq \max\{t>0\,;\,\Vert
b\Vert_{L^2(0,T_f)}\sqrt{t}\,(1+\frac{2\,R}{3}t\big)+\frac{2\sqrt{M\,t}}{\sqrt{\pi}}
+\frac{4\,R\,\sqrt{M}\,t^{3/2}}{3\,\sqrt\pi}
\,\leq 1\}\;.
\end{equation*}
Since we already have demanded that $t\leq T\leq \frac{1}{R}$ a 
sufficient condition is
then
$$\sqrt{T}\Big[\frac{5}{3}\,\Vert
b\Vert_{L^2(0,T_f)}+\frac{10\,\sqrt{1+2\alpha}}{3\,\sqrt{\pi}}
\Big]\,\leq 1\;,
$$
which reduces to  $T\leq T_c$ with $T_c$ given by \eqref{eq:def-tc}.\\

\medskip

\noindent
Our next step will consist in establishing  \textit{a priori} bounds 
on $p$ in $L^\infty_t(L^2_\sigma)\cap 
L^2_t(H^1_\sigma)$.\\

\vskip10pt
\noindent\textbf{Step 2}: \textit{\textsl{A priori} bounds.}

\medskip

\noindent
If we multiply equation (\ref{eq:plin}) by $p$ and
integrate by parts over $\RR$ with respect to $\sigma$ we easily 
obtain
\begin{eqnarray*}
\frac{1}{2}\,\frac{d}{dt}\int_{\RR}\,p^2\,+\,a\,\int_{\RR}|\partial_{\sigma}p|^2\,\le\,
\frac{D(q)}{\alpha}\,p(t,0)\;.
\end{eqnarray*}
Since from the Sobolev embedding of $H^1(\RR)$ into $L^\infty(\RR)$ 
and the bound
\eqref{eq:bd-D} on $D(q)$ we get
\begin{eqnarray*}
|\frac{D(q)}{\alpha}\,p(t,0)|&\le&R\,\Vert
 p\Vert_{L^\infty_\sigma}\\
&\le&R\,\big(\Vert
 p\Vert^2_{L^2_\sigma}\,+\,\Vert\partial_\sigma
 p\Vert^2_{L^2_\sigma}\big)^{\frac 1 2}\\
&\le&\frac{R^2}{2\,\varepsilon}\,+\,\frac \varepsilon 2
 \Vert p\Vert^2_{L^2_\sigma}\,+\,\frac \varepsilon 2
 \Vert \partial_\sigma p\Vert^2_{L^2_\sigma}\;,
\end{eqnarray*}
we may write
\begin{eqnarray}
\label{eq:endsigmap}
\frac{1}{2}\,\frac{d}{dt}\Vert p\Vert^2_{L^2_\sigma}\,+\,(a-\frac
\varepsilon 2)\,\Vert \partial_{\sigma}p\Vert^2_{L^2_\sigma}\,\le\,
\frac{R^2}{2\,\varepsilon}\,+\,\frac \varepsilon 2
 \Vert p\Vert^2_{L^2_\sigma}\;.
\end{eqnarray}
We recall that $a\geq \varepsilon$
and we apply the Gronwall lemma to obtain
\begin{equation}
\label{eq:bd-L2}
\sup_{0\leq t\leq T}\Vert 
p\Vert^2_{L^2_\sigma}\,\le\,e^{\varepsilon\,T}\big(\,\Vert
p_0\Vert^2_{L^2_\sigma}\,+\,\frac{T\,R^2}{\varepsilon}\,
\big)\;.
\end{equation}
We now return to  (\ref{eq:endsigmap}) and integrate it over $[0;T]$ 
to obtain
\begin{equation}
\varepsilon\,\Vert
\partial_{\sigma}p\Vert^2_{L^2_{t,\sigma}}\leq\Vert
p_0\Vert^2_{L^2_\sigma}\,(1\,+\,\varepsilon\,T  \,e^{\varepsilon
  T})\,+\frac{T\,R^2}{\varepsilon}\,(1\,+\,\varepsilon\,T  
\,e^{\varepsilon
  T}).\label{in:endsigmap}
\end{equation}

\vskip10pt
\noindent\textbf{Step 3}: \textit{The function ${\cal T}$ is 
continuous.}
\medskip

\noindent
We consider a sequence $(a_n;q_n)$ in ${\cal D}_{\varepsilon,M}\times 
Y_R$
such that $a_n$ converges
to  $a$ strongly in $L^2_{t}$
 and $q_n$ converges to $q$ strongly in $L^2_{t,\sigma}$, and   we 
denote
${\cal T}(a_n;q_n)\,=\,(D(p_n)+\varepsilon;p_n)$. We have to  prove 
that $p_n $ converges strongly
to $p$ in $L^2_{t,\sigma}$ and $D(p_n)$ converges to $D(p)$
strongly in $L^2_t$, with $(D(p)+\varepsilon\,; p)={\cal T}(a;q)$.

\medskip

\noindent
In virtue of (\ref{eq:bd-L2}) and (\ref{in:endsigmap}), the
sequence $p_n$ is bounded in
$L^\infty_t(L^{2}_\sigma)\,\cap\,L^2_t(H^1_\sigma)$. Then,
$\partial_\sigma p_n$ is bounded in $L^\infty_t(H^{-1}_\sigma)$
and $\partial^2_{\sigma\sigma}p_n$ is bounded in
$L^2_t(H^{-1}_\sigma)$. Since $a_n\partial^2_{\sigma\sigma}p_n$ is
bounded in $L^2_t(H^{-1}_\sigma)$, $b\,\in\,L^2_t$ and
$D(q_n)\delta_0$ is bounded in   $L^2_t(H^{-1}_\sigma)$,
$\partial_t p_n$ is bounded in $L^2_t(H^{-1}_\sigma)$. This
together with the fact that $p_{n}$ is bounded in
$L^2_t(H^1_\sigma)$ implies that, up to a subsequence, $p_n$
converges strongly towards $ p$  in $L^2_t(L^2_{{\rm
loc},\sigma})$ (the convergence being weak in $L^2_{t}(H^1_\sigma)$\/)
thanks to a well-known compactness result~\cite{bi:JLL}. In 
particular,
$p_{n}$  converges to $p$
almost everywhere. Thus $p\geq 0$ and  by the Fatou's lemma,
$\int_{\RR}|\sigma|\,p\leq R$ almost everywhere on $[0;T]$. Hence
$p$ belongs to $Y_R$. We  are going to show that the
convergence is actually strong in $L^2_{t,\sigma}$.

\medskip

\noindent
In virtue of
\eqref{eq:bd-infini} in Proposition~\ref{prop:1},  we dispose of a
uniform \textsl{a 
priori\/} bound on
$p_n$ in
$L^\infty_{t,\sigma}$ (hence also on $p$). For the strong convergence 
in $L^2_{t,\sigma}$ we then
argue as follows. For any fixed positive real number
$K$, we have
\begin{eqnarray*}
\int_0^T\int_\RR \vert 
p_n-p\vert^2&\le&\int_0^T\int_{\vert\sigma\vert\le K}\vert
p_n-p\vert^2\,+\,\int_0^T\int_{\vert\sigma\vert>K}\vert p_n-p\vert^2\\
&\le&\int_0^T\int_{\vert\sigma\vert\le K}\vert
p_n-p\vert^2\,+\,\big(\Vert
p_n\Vert_{L^\infty_{t,\sigma}}+\Vert
p\Vert_{L^\infty_{t,\sigma}}\big)\frac{2\,R\,T}{K}\;,
\end{eqnarray*}
owing to the fact that $p_n$ and $p$ belong to a bounded subset of
$Y_R\cap L^\infty_{t,\sigma}$. We then  conclude  by letting
$n$ next $K$ go to infinity.\\

\medskip

\noindent
We now prove that
$D(p_n)$ converges to
$D(p)$ strongly in
$L^2_{t}$. We shall actually prove that $D(p_n)$ converges to
$D(p)$ strongly in
$L^1_{t}$ and then  use the  fact that $D(p_{n})$ is bounded in 
$L^\infty_t$, in virtue
of
\eqref{eq:bd-D} and because  $p_{n}$ lies in $Y_R$. Let us fix $K>1$. 
Then,  we have
\begin{eqnarray}
\frac 1\alpha \int_0^T
|D(p_n)-D(p)|&=&\int_0^T\big|\int_{|\sigma|>1}(p_n\,-\,p)\,\big|\nonumber\\
&\le&\int_0^T\int_{1<|\sigma|<K}|p_n\,-\,p|\,
+\,\frac{1}{K}\,\int_0^T\int_{|\sigma|>K}|\sigma|\,\big(|
p_n|+|p|\big)\nonumber\\
&\le&\int_0^T\int_{1<|\sigma|<K}|p_n\,-\,p|\,+\frac{2\,R\,T}{K}\;,\label{in:convgeD}
\end{eqnarray}
because $p$ and $p_{n}$ belong to $Y_R$. Since $p_n$ converges  to 
$p$ strongly in
$L^1_t(L^1_{{\rm loc},\sigma})$, we conclude that $D(p_n)$ converges 
to $D(p)$ in  $L^1_t$ by
letting $n$ next $K$ go to infinity in (\ref{in:convgeD}).\\

\medskip

\noindent
 In order to pass to the limit in the equation satisfied by $p_n$ 
(thereby
proving that $(D(p)+\varepsilon\,; p)={\cal T}(a;q)$\/), we now
observe that the strong convergence of $q_n$ to $q$ in 
$L^2_{t,\sigma}$, together with
the argument in \eqref{in:convgeD} above shows that
$D(q_{n})$ converges to $D(q)$ strongly in $L^{2}_t$. It is then 
easily proved that $p$ is a weak solution 
to~\eqref{eq:p-eps} and since $p$ is in $L^2_t(H^1_\sigma)$ it  is the
unique solution to \eqref{eq:syst-lin} corresponding to
$a$ and $q$. In particular, the whole sequence $p_{n}$ converges and
not only a subsequence.
\vskip10pt
\noindent\textbf{Step 4}: \textit{${\cal T}({\cal 
D}_{\varepsilon}\times
  Y_R)$ is relatively compact.}\\

\medskip

\noindent
Let $(D(p_n)+\varepsilon\,; p_{n})={\cal T}(a_{n};,q_{n})$ be a
sequence in ${\cal T}({\cal D}_{\varepsilon,M}\times Y_R)$. We
have to prove that we may extract a subsequence which converges
strongly in $L^{2}_t\times L^{2}_{t\,,\sigma}$. Exactly as for the
proof of the continuity, the \textsl{a priori} estimates
(\ref{eq:bd-L2}) and (\ref{in:endsigmap})  ensure that the
sequence $p_n$ is bounded in
$L^\infty_t(L^{2}_\sigma)\,\cap\,L^2_t(H^1_\sigma)$. Since
$|\sigma|\, p_{n}$ is bounded $L^{\infty}_t(L^{1}_\sigma)$, we can
mimic the argument in Step~3  above to deduce that up to a
subsequence the sequence  $p_{n}$ converges to some $p$ in $Y_R$
strongly in $L^{2}_{t\,,\sigma}$ and  that $D(p_n)$ converges to $
D(p)$ strongly in $L^2_{t}$. \hfill $\diamondsuit$ 

\medskip

\noindent
We are now in position to conclude the proof of 
Proposition~\ref{prop:visc}. \\

\medskip

\noindent
Let $T_f>0$ and $0<\varepsilon\leq 1$ being given. We are going to 
prove the existence of a unique
solution on $[0;T_f]$.\\

\medskip

\noindent
Being given an initial data $p_0$ which satisfies \eqref{eq:IC-full},
existence  of a solution 
$p_\varepsilon$ is ensured  from Proposition~\ref{prop:2} by applying
the Schauder fixed point theorem 
on ``short'' time interval $[0;T_1]$ with 
$T_1=\min(\frac{\alpha}{R_1},
T_c)$ and where $R_1=1+\int_\RR 
\vert\sigma\vert p_0\,$. This solution is uniquely defined in virtue 
of
Lemma~\ref{lem:unique} and 
we know from \eqref{eq:proba} that $\int_\RR
p_{\varepsilon}(T_1)=1$. Moreover from 
Proposition~\ref{prop:1} $p_\varepsilon(T_1)\in L^\infty_{\sigma}$ and
by construction 
$\int_\RR \,|\sigma|\,p_\varepsilon(T_1)\leq R_1$. Therefore
$p_\varepsilon(T_1)$ satisfies the same conditions \eqref{eq:IC} as
$p_0$.  Then, repeating the same argument we may build a solution to
\eqref{syst:p-eps} with initial data $p_\varepsilon(T_1)$ on 
$[T_1;T_2]$ with
$T_2=\min\big(\frac{1}{R_2},T_c\big)$, where $R_2=R_1+1=\int_\RR
|\sigma|\,p_0+2$. Thanks to the 
uniqueness result (Lemma~\ref{lem:unique}), if we now glue this 
solution
to $p_\varepsilon$ at $t=T_1$ we 
obtain the unique solution to \eqref{syst:p-eps} on $[0;T_1+T_2]$. It 
is
now clearly seen that for any 
integer $n\geq 1$ we may build a solution to \eqref{syst:p-eps} on
$[0;\sum_{1\leq k\leq n}T_k]$ with 
$T_k=\min\Big(\big(k+\int_\RR
|\sigma|\,p_0\big)^{-1};T_c\Big)$. Since $\sum_{1\leq k\leq 
n}T_k$ obviously goes to $+\infty$ together with $n$, existence (and
uniqueness) of the  solution 
$p_\varepsilon$  to
\eqref{syst:p-eps} is obtained  on every time interval.\\

For the proof of \eqref{eq:bd-infini-eps} we argue as for the proof 
of  \eqref{eq:bd-infini}  in Proposition~\ref{prop:1}.
Defining $p_\varepsilon^+$ as in \eqref{eq:p+} with $a$ replaced by 
$D(p_\varepsilon)+\varepsilon$ and
$D(q)$ by $D(p_\varepsilon)$ we obtain
\begin{eqnarray*}
0&\leq &p_\varepsilon\leq p_\varepsilon^+\\
&\leq &\Vert p_0\Vert_{L^\infty}+\frac{1}{\alpha\,\sqrt\pi}\int_0^t 
\frac{D\big(p_\varepsilon(s)\big)}
{2\,\sqrt{\varepsilon+\int_s^tD(p_\varepsilon)}}\,ds\\
&\leq &\Vert 
p_0\Vert_{L^\infty}+\frac{1}{\alpha\,\sqrt\pi}\;\Big[\sqrt{\varepsilon+\int_0^tD(p_\varepsilon)}-\sqrt\varepsilon\Big]\\
&\leq &\Vert 
p_0\Vert_{L^\infty}+\frac{1}{\alpha\,\sqrt\pi}\;\sqrt{\int_0^t 
D(p_\varepsilon)}\\
&\leq&\Vert p_0\Vert_{L^\infty}+ \frac{\sqrt\alpha\,\sqrt 
T}{\sqrt{\pi}}\;.
\end{eqnarray*}
Then
\begin{equation*}
\int_\RR p_\varepsilon^2\leq \Vert 
p_\varepsilon\Vert_{L^\infty_\sigma}\;\int_\RR
p_\varepsilon\;,\end{equation*}
from which \eqref{eq:bd-L2-eps} follows gathering together 
\eqref{eq:proba} and
\eqref{eq:bd-infini-eps} and,  with the notation of the proposition, 
\begin{equation*}
C_2(T,p_0)=\Vert p_0\Vert_{L^\infty}+ \frac{\sqrt\alpha\,\sqrt 
T}{\sqrt{\pi}}\;.
\end{equation*} \\

\medskip

\noindent
The proof of \eqref{eq:bd-sigma-eps} follows the same lines as the 
proof of \eqref{eq:bd-sigma}.
Indeed, we again use the pointwise inequality $p_\varepsilon\leq 
p_\varepsilon^+$  and replace $D(q)$
by $D(p_\varepsilon)\,(\leq \alpha)$ and $a$ by 
$D(p_\varepsilon)+\varepsilon\,( \leq
\alpha+1)$ in  \eqref{eq:bd-sigma} and use \eqref{eq:IC} to deduce
\begin{equation}\label{eq:etape}
\sup_{0\leq t\leq T}\int_\RR\vert\sigma\vert\, p_\varepsilon\leq
\int_\RR \vert\sigma\vert\, 
p_0+\sqrt{T}\,\Big(\frac{2\,\sqrt{1+\alpha}}{\sqrt{\pi}}+\Vert
b\Vert_{L^2(0,T)}\Big)+\frac{2}{3}\,T^{3/2}\,\Big(1+\frac{2\sqrt{1+\alpha}}{\sqrt\pi}\Big)\;,
\end{equation}
whence \eqref{eq:bd-sigma-eps} with $C_1(T,p_0)$ being the
quantity in the right-hand side of \eqref{eq:etape}.

\medskip

\noindent
In order to prove \eqref{eq:bd-grad-eps},  we  multiply  
(\ref{eq:p-eps})
by  $p_\varepsilon,$ and we integrate by parts over
$\RR$  with respect to
$\sigma$ to obtain
\begin{equation}\label{*}
\frac{1}{2}\frac{d}{dt}\int_{\RR}  p_\varepsilon^2
\,+\,\big(D(p_\varepsilon)+\varepsilon\big)\,\int_{\RR} 
\vert\partial_{\sigma} p_\varepsilon
\vert^2
\,+ \int_{\vert \sigma \vert>1} p_\varepsilon^2
\,=\,\frac{D(p_\varepsilon)}{\alpha}\; p_\varepsilon(t,0)\;.
\end{equation}
We use the $L^\infty$ bound \eqref{eq:bd-infini-eps} to bound the 
right-hand side and we integrate
\eqref{*} with respect to $t $ over $[0;T]$ to deduce 
\eqref{eq:bd-grad-eps} with
\begin{equation*}
C_3(T,p_0)=\Vert p_0\Vert_{L^\infty}\big(\frac 1 2+ 
T)+\,\frac{\sqrt{\alpha}}{\sqrt\pi}\,T^{3/2}\;,
\end{equation*}
using that $\Vert
p_0\Vert_{L^2_{\sigma}}^2\leq \,\Vert
p_0\Vert_{L^\infty}\int_\RR p_0$.
\hfill $\diamondsuit$
%
%
\section{The non degenerate case: $D(p_0)>0$}\label{sec:nondeg}

The main result of this section  corresponds to the statement of 
Theorem~\ref{th:main1}  and
fully describes the issue of existence and uniqueness of solutions to
the HL equation 
\eqref{eq:syst-p} in the non-degenerate case. It is summarized in  the
following~: 
\begin{proposition}\label{prop:non-deg} Let $p_0 $ satisfy
\eqref{eq:IC-full}. We assume that  $D(p_0)>0$.  Then, the HL equation
\eqref{eq:syst-p} has a unique solution
  $p$ in $C^0_t(L^2_\sigma)\cap L^2_t(H^1_\sigma)$ and  $p$ is
the  limit (in $L^2_{t,\rm{loc}}(L^2_\sigma) \cap 
C_{t,\rm{loc}}^0(L^2_\sigma)$) of
$(p_\epsilon)$ when $\epsilon$ 
  goes to $0$ where $p_\varepsilon$ is the vanishing viscosity 
solution
  whose existence and uniqueness 
is ensured by Proposition~\ref{prop:visc}. Moreover, $p\in
L^\infty_{t,\sigma}\cap C^0_t(L^1_\sigma)$, 
$\sigma\,p\in L^\infty_t(L^1_\sigma)$ and $\int_\RR p=1$. Furthermore,
$D(p)\in C^0_t$ and for every 
$T>0$ there exists a positive constant $\nu(T)$ such that
\begin{equation}\label{eq:strict-pos}
\min_{0\leq t\leq T}D(p(t))\geq \nu(T)\;.
\end{equation}
\end{proposition}
\medskip

\noindent
We begin with proving the following~:

\begin{lemme}\label{lem:non-degeneracy} We assume that $p_0$ 
satisfies \eqref{eq:IC-full}.
Then, if $D(p_0)>0$, $D(p_\varepsilon)(t)>0$ for every $t\in [0,T]$, 
with $p_\varepsilon$ being the
unique solution to \eqref{syst:p-eps} provided by 
Proposition~\ref{prop:visc} and,  actually,
for every
$T>0$ there exists a positive constant $\nu(T)$ such that
\begin{equation}\label{eq:strict-pos-eps}
\min_{0\leq t\leq T}D(p_\varepsilon(t))\geq \nu(T)\;,
\end{equation}
for every $0< \varepsilon\leq 1$. 
\end{lemme}

\begin{remarque} Note that this bound from below is independent of 
$\varepsilon$, but it
comes out from the proof that it depends  on $p_0$ and on  the shear  
$b$.
\end{remarque} 
\noindent\textbf{Proof of Lemma~\ref{lem:non-degeneracy}:} The proof 
relies on the bound from
below in \eqref{eq:compare} that we integrate over $\vert \sigma 
\vert>1$ to obtain
\begin{equation}
D(p_\varepsilon(t))\geq\alpha\int_{\vert\sigma\vert>1} p_\varepsilon^-
\geq
\alpha\, e^{
-t }\,
\int_\RR p_{0}(\sigma')\,\Big( \int_{\vert \sigma
\vert>1}
\varphi_{\sqrt{2\int_0^t 
(D(p_\varepsilon)+\varepsilon)}}(\sigma-\sigma'-\chi(t))
 \,d\sigma\Big)\,d\sigma'\;.\label{eq:strict}
\end{equation}
Let us define
$\displaystyle{K_\chi\,=\,[\,-1-\chi(t),1-\chi(t)\,]}$.
The function $\sigma\mapsto \varphi_{\sqrt{2\int_0^t
(D(p_\varepsilon)+\varepsilon) }}(\sigma-\sigma'-\chi(t))$ is a 
Gaussian probability density with mean
$\sigma'+\chi(t)$ and
squared width $2\,\int_{0}^{t}\,(D(p_\varepsilon)+\varepsilon)$. 
Therefore, for every
$\displaystyle{\sigma'\in\RR\setminus K_{\chi}}$,  we have
$$\int_{\vert\sigma \vert>1}\,\varphi_{\sqrt{2\int_0^t
(D(p_\varepsilon)+\varepsilon)}}(\sigma-\sigma'-\chi(t))\,d\sigma\ge\frac{1}{2}\;,$$

which implies
\begin{equation*}
\eqref{eq:strict}\geq \frac{\alpha}{2}\,e^{-T}\,\int_{\RR\setminus 
K_{\chi}}\,p_0=
\frac{\alpha}{2}\,e^{-T}
 \,\int_{\vert\sigma+\chi(t)\vert>1} p_0\;.
\end{equation*}
In the zero shear case ($b\equiv 0$, thus
$\chi\equiv 0$) the proof is over and  
\begin{equation*}
\min_{0\leq t\leq T}D(p(t))\geq \frac{1}{2}e^{-T}D(p_0)\;.
\end{equation*}
In the general case, a strictly positive bound from below  is 
available as long as 
the support of $p_0$ is not contained in $K_\chi$. We thus define 
\begin{equation}\label{def:t*}
t^* =\inf\Big\{t>0\,; \int_{\vert\sigma+\chi(t)\vert>1} 
p_0=0\,\Big\}\;.
\end{equation}
Then $0<t^*$ ($t^*$ possibly even  infinite),  the support of $p_0$
is contained in 
$[-1-\chi(t^*),1-\chi(t^*)[$, and
for every
$T<\frac{t^*}{2}$,
\eqref{eq:strict-pos-eps} holds for some positive constant $\nu_1(T)$ 
defined by
\begin{equation}\label{eq:nu1}
\nu_1(T)= \frac{\alpha}{2}\,e^{-T}\,\min_{0\leq t\leq T}
 \,\int_{\vert\sigma+\chi(t)\vert>1} p_0.
\end{equation}
It is worth emphasizing that this quantity is independent of 
$\varepsilon$. If
$t^*=+\infty$, the proof is over and $\nu(T)=\nu_1(T)$ fits. Let us 
now
examine the case 
when $t^*<+\infty$ and $T\geq \frac{t^*}{2}$.  \\

\noindent
We go back to \eqref{eq:strict}, take $t$ in $[\frac{t^*}{2};T]$ and 
denote $x=\int_0^t
(D(p_\varepsilon)+\varepsilon)$ for shortness. Then
\begin{eqnarray}
 D(p_\varepsilon(t))& \geq &\alpha\, e^{
-T }\,
\int_{-1-\chi(t^*)}^{1-\chi(t^*)}p_{0}(\sigma')\,\Big( \int_{\vert 
\sigma
\vert>1}
\varphi_{\sqrt{2\,x}}(\sigma-\sigma'-\chi(t))
 \,d\sigma\Big)\,d\sigma'\nonumber\\
& = & \alpha\, e^{
-T }\,\int_{-1-\chi(t^*)}^{1-\chi(t^*)} p_0(\sigma')
 \left(  \int_{|\sigma|>1} 
\frac{e^{-(\sigma-\sigma'-\chi(t))^2/4x}}{2 \, \sqrt{\pi}
     \, \sqrt x} d\sigma \right) \, d\sigma' \nonumber\\
& = & \frac {\alpha} {\sqrt \pi}e^{-T }\,
\int_{-1-\chi(t^*)}^{1-\chi(t^*)} 
p_0(\sigma')\left(\int_{-\infty}^{-1+\sigma'+\chi(t)}
  \frac{e^{-\sigma^2/4x}}{2 \, \sqrt x} d\sigma +
\int_{1+\sigma'+\chi(t)}^{+\infty} \frac{e^{-\sigma^2/4x}}{2 \,  
\sqrt x}
d\sigma \right) \,
d\sigma' \nonumber\\ 
& = &  \frac {\alpha} {\sqrt \pi}e^{-T }\,
\int_{-1-\chi(t^*)}^{1-\chi(t^*)} p_0(\sigma')\,\left(
\int_{\frac{1+\sigma'+\chi(t)}{2\sqrt x}}^{+\infty} e^{-t^2} \, dt +
\int_{\frac{1-\sigma'-\chi(t)}{2\sqrt x}}^{+\infty} e^{-t^2} \, dt 
\right) \,
d\sigma'\nonumber\\
&\geq& \frac{\alpha}{\sqrt \pi}e^{-T }\,\left(
\int_{-1-\chi(t^*)}^{1-\chi(t^*)}
p_0(\sigma')\,d\sigma'\right)\,\left(\int_{\frac{2-\chi(t^*)+\chi(t)}{\sqrt{2\,t^*\,\nu_1(t^*/2)}}}^{+\infty}
e^{-t^2}\, dt +
\int_{\frac{2+\chi(t^*)-\chi(t)}{\sqrt{2\,t^*\,\nu_1(t^*/2)}}}^{+\infty} 
e^{-t^2} \,
dt\right)\nonumber\\
&\geq& \frac{\alpha}{\sqrt \pi}e^{-T
}\min_{t^*/2\leq t\leq
T}\,\left(\int_{\frac{2-\chi(t^*)+\chi(t)}{\sqrt{2\,t^*\,\nu_1(t^*/2)}}}^{+\infty} 
e^{-t^2}\, dt +
\int_{\frac{2+\chi(t^*)-\chi(t)}{\sqrt{2\,t^*\,\nu_1(t^*/2)}}}^{+\infty} 
e^{-t^2} \,
dt\right)\;,\label{eq:ouf}
\end{eqnarray}
since $ \int_{-1-\chi(t^*)}^{1-\chi(t^*)}
p_0=1$ and  $x\geq \int_0^{t^*/2}D(p_\varepsilon)\geq 
t^*\,\nu_1(t^*/2)/2$
thanks to
\eqref{eq:nu1}. The proof of Lemma~\ref{lem:non-degeneracy} then 
follows by
defining
\begin{equation*}
\nu(T)=\min(\nu_1(T);\nu_2(T))\;,
\end{equation*}
with $\nu_1(T)$ given by \eqref{eq:nu1} and $\nu_2(T) $ being the 
positive quantity in the
right-hand side of
\eqref{eq:ouf}, that is 
\begin{equation*}
\nu_2(T)=\frac{\alpha}{\sqrt \pi}e^{-T
}\min_{t^*/2\leq t\leq
T}\,\left(\int_{\frac{2-\chi(t^*)+\chi(t)}{\sqrt{2\,t^*\,\nu_1(t^*/2)}}}^{+\infty} 
e^{-t^2}\, dt +
\int_{\frac{2+\chi(t^*)-\chi(t)}{2\sqrt{2\,t^*\,\nu_1(t^*/2)}}}^{+\infty} 
e^{-t^2} \,dt\right)\;.
\end{equation*}
{}\hfill$\diamondsuit$
\vskip10pt\noindent\textbf{Proof of Proposition~\ref{prop:non-deg}:}\\

\medskip

\noindent
We first  go back to the proof of the bound \eqref{eq:bd-grad-eps} on
$\partial_\sigma p_\varepsilon$ and more precisely we  look at
\eqref{*} and observe that in virtue of \eqref{eq:strict-pos}
\begin{equation}\label{eq:bd-grad}
\nu(T)\,\int_0^T\int_\RR\vert\partial_\sigma\,p_\varepsilon\vert^2\leq 
C_3(T,p_0)\;.
\end{equation}
Let now $\varepsilon_n$ denote any sequence in $[0,1]$ which goes to 
$0$ as $n$ goes to infinity. To
shorten the notation we denote by $p_n$ instead of $p_{\varepsilon_n}$
the corresponding sequence 
of solutions to \eqref{syst:p-eps}. With the above bound 
\eqref{eq:bd-grad} on $p_n$ and
\eqref{eq:bd-L2-eps}, we know that $p_n$ is bounded in 
$L^2_t(H^1_\sigma)$ independently of $n$.
Moreover thanks to
\eqref{eq:proba} and \eqref{eq:bd-infini-eps}
$p_n$ is bounded  in $L^\infty_t(L^1_\sigma\cap L^\infty_\sigma)$ and 
we also dispose of a uniform
bound  on $\int_\RR \vert\sigma\vert p_n$ in virtue of 
\eqref{eq:bd-sigma-eps}. Therefore arguing
exactly as in the proof of Proposition~\ref{prop:2} (Step 4) where we 
have proved that the mapping ${\cal T}
$ is  relatively compact in $L^2_t\times L^2_{t,\sigma}$ we show that 
$p_n$ converges to some $p$ strongly 
in $L^2_{t,\sigma}$ and $D(p_n)$ converges to $D(p)$ in $L^2_t$. Then
$p$ is a solution to the initial problem 
\eqref{eq:syst-p}  in $L^2_t(H^1_\sigma)\cap 
L^\infty_t(L^1_\sigma\cap L^\infty_\sigma)$, $\int_\RR p=1$ and
$\int_\RR\vert\sigma\vert\,p<+\infty$. Moreover, 
\begin{equation*}
\inf_{0 \le t \le T} D(p(t))\geq \nu(T)\;.
\end{equation*}
This non-degeneracy condition on the viscosity coefficient ensures 
that there is at most
one solution to \eqref{eq:syst-p} in $L^2_t(H^1_\sigma)\cap 
L^\infty_t(L^2_\sigma)$ (this follows by  an
obvious adaptation of the proof  of Lemma~\ref{lem:unique} to this 
case). Therefore the limiting
function
$p$ is uniquely defined and does not depend on the sequence 
$\varepsilon_n$. Moreover the whole
sequence
$p_n$ converges to this unique limit and not only a
subsequence. \hfill$\diamondsuit$ \\ 

\medskip

\noindent
As a conclusion of this subsection let us make the following comment 
which is a byproduct of
Proposition~\ref{prop:non-deg}. Let $p$ be a solution to 
\eqref{eq:syst-p} in $C^0_t(L^1_\sigma\cap L^2_\sigma)$, then as soon 
as $D(p(t))$ is
positive for some time $t$ it remains so afterwards since the 
solution can be continued 
in a unique way beginning from time $t$. 

\vskip10pt
%
%
\section{The degenerate case~: $D(p_0)=0$}\label{sec:deg}
\medskip
Throughout this section we assume that $D(p_0)=0$ and therefore the
support of $p_0$ is 
included in $[-1;+1]$. Assume that we dispose of a solution
to~\eqref{eq:syst-p} in 
$C^0_t(L^1_\sigma\cap L^2_\sigma)$. We may define
$t_*\in
\RR^+ \cup\{+\infty\}$ by
\begin{equation}\label{eq:def-t*}
t_*=\max\Big\{t>0\, ; \int_0^tD(p )=0\Big\}\;.
\end{equation}
According to the comment at the end of the previous section for every 
$t>t_*$, $D(p(t))>0$ while $D(p(t))=0$ for all $t$ in $[0;t_*]$. On $[0;t_*[$,
the HL equation~\eqref{eq:syst-p} reads
\begin{equation*}
\left\{
\begin{array}{rcl}
\partial_t\,p&=&-b(t)\,\partial_\sigma p\;;\\
p&\geq& 0\;;\\
p(0,\cdot)&=&p_0\;;\\
D(p(t))&=&0\;.
\end{array}
\right.
\end{equation*}
The above system reduces to 
\begin{equation}\label{eq:transp}
\left\{
\begin{array}{rcl}
p(t,\sigma)&=&p_0\big(\sigma-\chi(t)\big)\;;\\
D(p(t))&=&0\,,\quad \textrm{for all $t$ in $[0;t_*]$ }\;.
\end{array}
\right.
\end{equation}
The second equation in \eqref{eq:transp} is compatible with the first
one as long as  
\begin{equation*}
\int_{\vert \sigma+\chi (t)\vert>1} p_0=0\,,\quad \textrm{for all $t$ 
in $[0;t_*]$ }\;. 
\end{equation*}
Therefore there exists a maximal time interval $[0;T_c]$ on which the 
HL equation may
reduce to a mere transport equation and this is for an intrinsic time 
$T_c$ (possibly
infinite) defined by 
\begin{equation}\label{eq:defTc}
T_c=\inf \left\{ t>0\,;\; \int_{\vert \sigma+\chi (t)\vert>1} p_0>0 
\;.
\right\}\;.  
\end{equation}
Note that $T_c$ is completely determined by the data $p_0$ and $b$. If
$T_c=+\infty$, the 
steady state $p(t,\sigma)=p_0(\sigma-\chi(t))$ is a solution of the HL
equation for all time. We shall now exhibit circumstances under which 
it
is not the unique solution. For convenience, we restrict ourselves to 
the
case when $b \equiv 0$ (we then have obviously $T_c=+\infty$).


\medskip

\noindent
For $p_0 \in L^1 (\RR) \cap L^\infty(\RR)$ such that $p_0 \ge 0$, let 
us denote by $F_{p_0}$ the function from $\RR^+$ to $\RR^+$ defined by
$F_{p_0}(0) = D(p_0)$ and by
\begin{equation*}
\forall x > 0 , \qquad
F_{p_0}(x) = \alpha \, \int_{|\sigma|>1} \left( \int_{\RR} 
p_0(\sigma')
\varphi_{\sqrt{2x}}(\sigma-\sigma') \, d\sigma' \right) \,d\sigma\;.
\end{equation*}

\medskip

\noindent
\begin{proposition}\label{prop:deg} Let $p_0 $ satisfy 
\eqref{eq:IC-full} and 
be  such that $D(p_0) = 0$, then
\begin{enumerate}
\item If $F_{p_0}$ satisfies
\begin{equation}
 \label{eq:condFp0}
\int_0^1 \frac{dx}{F_{p_0}(x)} = + \infty,
\end{equation}
then $p(t,\sigma)=p_0(\sigma)$ is the unique solution to
\eqref{eq:syst-p} in $C^0_t(L^2_\sigma)$~;
\item Otherwise, \eqref{eq:syst-p} has an infinite number of 
solutions in
  $C^0_t(L^2_\sigma)$. The set of solutions to \eqref{eq:syst-p} is 
made of the
  steady state $p(t,\sigma)=p_0(\sigma)$ and of the functions
$\left( q_{t_0} \right)_{t_0 \ge 0}$ defined by
\begin{equation*}
q_{t_0}(t,\sigma) \, = \, \left| \begin{array}{lll}
\dps p_0(\sigma) & \quad & \mbox{if } t \le t_0 \\
\dps q(t-t_0,\sigma) & & \mbox{if } t > t_0 \end{array} \right.
\end{equation*}
where $q$ is the unique solution to \eqref{eq:syst-p} in
$C^0_t(L^2_\sigma)$ such that 
$D(q)>0$ on $]0,+\infty[$. Besides,
\begin{equation}
p_\epsilon \mathop{\longrightarrow}_{\epsilon \rightarrow 0} q \qquad
\mbox{strongly in $L^2_{t,\rm{loc}}(L^2_\sigma)$}.
\end{equation}
\end{enumerate}
\end{proposition}
\medskip

\begin{lemme}\label{lem:1} Let $p_0 \in L^1 (\RR) \cap L^\infty(\RR)$ 
such that
\begin{equation*}
p_0 \ge 0, \qquad \int_\RR p_0 = 1, \qquad D(p_0)=0.
\end{equation*}
The function $F_0$ is in $C^0([0,+\infty[) \cap 
C^\infty(]0,+\infty[)$,
and is positive on $]0,+\infty[$. In addition, $F_{p_0}'> 0$ on 
$]0,+\infty[$.
\end{lemme}
\medskip

\noindent\textbf{Proof of Lemma~\ref{lem:1}:} It is easy to check
that $F_{p_0} \in C^0([0,+\infty[) \cap C^\infty(]0,+\infty[)$,
and that $F_{p_0} > 0$ on $]0,+\infty[$. Since $D(p_0)=0$, the 
function
$p_0$ is supported 
in
$[-1,1]$. Thus, for any $x > 0$
\begin{eqnarray}
F_{p_0}(x) & = & \alpha\,\int_{|\sigma|>1} \left( \int_\RR 
p_0(\sigma')
  \varphi_{\sqrt{2x}}(\sigma-\sigma') \, d\sigma' \right) \, d\sigma 
\nonumber\\
& = & \alpha\int_{-1}^1  p_0(\sigma')
 \left(  \int_{|\sigma|>1} \frac{e^{-(\sigma-\sigma')^2/4x}}{2 \, 
\sqrt{\pi}
     \, \sqrt x} d\sigma \right) \, d\sigma' \nonumber\\
& = & \alpha\int_{-1}^1  p_0(\sigma') 
\left(\int_{-\infty}^{-1+\sigma'}
  \frac{e^{-\sigma^2/4x}}{2 \, \sqrt{\pi} \, \sqrt x} d\sigma +
\int_{1+\sigma'}^{+\infty} \frac{e^{-\sigma^2/4x}}{2 \, \sqrt{\pi} \, 
\sqrt x}
d\sigma \right) \,
d\sigma'\nonumber \\ & = & \alpha \frac 1 {\sqrt \pi} \int_{-1}^1  
p_0(\sigma') \left(
\int_{\frac{1+\sigma'}{2\sqrt x}}^{+\infty} e^{-t^2} \, dt +
\int_{\frac{1-\sigma'}{2\sqrt x}}^{+\infty} e^{-t^2} \, dt \right) \,
d\sigma'\;.\label{eq:defF}
\end{eqnarray}
It follows that for any $x > 0$,
\begin{equation*}
F_{p_0}'(x) = \alpha \frac 1 {\sqrt \pi} \int_{-1}^1  p_0(\sigma')
\left( \frac{1+\sigma'}{4x^{3/2}} e^{-\frac{(1+\sigma')^2}{4x}} +
\frac{1-\sigma'}{4x^{3/2}} e^{-\frac{(1-\sigma')^2}{4x}}
\right) \, d\sigma' > 0\;.
\end{equation*}
\hfill$\diamondsuit$

\medskip

\begin{lemme}\label{lem:2} Let $\gamma \geq 0$ and $p_0 \in L^1(\RR) 
\cap L^\infty(\RR)$
such that
\begin{equation*}
p_0 \ge 0, \qquad \int_\RR p_0 = 1, \qquad \int_\RR |\sigma| p_0 < +
\infty, \qquad D(p_0)=0.
\end{equation*}
Let us consider the problem
\begin{equation} \label{eq:II}
\left\{ \begin{array}{l}
\dps \partial_t w = D(w(t)) \partial^2_{\sigma \sigma} w -  \gamma w  
\\
w(0,\sigma) = p_0(\sigma). \end{array} \right.
\end{equation}
\begin{enumerate}
\item If $F_{p_0}$ satisfies (\ref{eq:condFp0}) then
  $p(t,\sigma)=p_0(\sigma)$ is the unique solution to (\ref{eq:II}) in
  $C^0_t(L^2_\sigma)$;
\item Otherwise, (\ref{eq:II}) has an infinite number of solutions in
  $C^0_t(L^2_\sigma)$. The set of solutions to (\ref{eq:II}) is made 
of the
  steady state $w(t,\sigma)=p_0(\sigma)$ and of the functions $\left(
    v_{t_0} \right)_{t_0 \ge 0}$ defined by
\begin{equation*}
v_{t_0}(t,\sigma) \, = \, \left| \begin{array}{lll}
\dps p_0(\sigma) & \quad & \mbox{if } t \le t_0 \\
\dps v(t-t_0,\sigma) & & \mbox{if } t > t_0 \end{array} \right.
\end{equation*}
where $v$ is the unique solution to
\eqref{eq:II} in $C^0_t(L^2_\sigma)$ such that $D(v)>0$ on 
$]0,+\infty[$.
\end{enumerate}
\end{lemme}
\medskip
\begin{corollaire}\label{cor:1} The initial data 
$p_0=\frac{1}{2}\chi_{]-1,1[}$ fulfills the assumptions of the above
lemma and $\int_0^1\frac{dx}{F_{p_0}(x)}<+\infty$. Therefore there 
are infinitely many solutions to the
equation~\eqref{eq:simpl}  in the introduction. 
\end{corollaire}
\noindent\textbf{Proof of Corollary~\ref{cor:1}: }
The only point to be checked is that 
$\int_0^1\frac{dx}{F_{p_0}(x)}<+\infty$. With the standard notation 
$\erfc(z)\equiv \int_z^{+\infty} e^{-t^2}\,dt$, and by using 
\eqref{eq:defF} and symmetry considerations, simple  calculations 
yield
\begin{eqnarray*}
F_{p_0}(x)&=&\frac{2\,\alpha\sqrt{x}}{\sqrt{\pi}}\,\int_0^{\frac{1}{\sqrt{x}}} 
\erfc(\sigma)\,d\sigma\\
&=& \frac{2\,\alpha}{\sqrt{\pi}}\,\left[\erfc\big(\frac 1 
{\sqrt{x}}\big)-\frac 1 2 \,\sqrt{x}\,e^{-\frac{1}{x}}+ \frac 12 
\,\sqrt{x}\right].
\end{eqnarray*}
Since $\erfc(z)\sim\frac{1}{2}e^{-z^2}/z $ for $z$ going to  
$+\infty$, $F_{p_0}(x)\sim \frac{\alpha}{\sqrt{\pi}} \,\sqrt{x}$ near 
$0$ and the integrability  of 
$1/F_{p_0}$ on $[0;1]$  follows. 
\hfill$\diamondsuit$

\bigskip
\noindent\textbf{Proof of Lemma~\ref{lem:2}}

\medskip

\noindent
Let us consider a non-negative function $D \in 
C^0([0,+\infty[,\RR^+)$.
The unique solution in $C^0_t(L^2_\sigma)$ of the problem
\begin{equation} \label{eq:III}
\left\{ \begin{array}{l}
\dps \partial_t w_D = D(t) \partial^2_{\sigma \sigma} w_D -  \gamma 
w_D; \\
w_D(0,\sigma) = p_0(\sigma), \end{array} \right.
\end{equation}
is given by
\begin{equation} \label{eq:wD}
w_D(t,\sigma) = \left| \begin{array}{lll}
e^{-\gamma t} \, p_0(\sigma) & \quad & \mbox{if } t \le t^\ast \, ; \\
&&\\ 
\dps e^{-\gamma t} \, \int_\RR p_0(\sigma') \varphi_{\sqrt{2
    \int_0^t D(s) \, ds}}(\sigma-\sigma') \, d\sigma' & \quad & 
\mbox{if } t >
t^\ast \, , \end{array} \right.
\end{equation}
where $t^\ast = \inf \left\{ t > 0, \quad \dps \int_0^t D > 0 
\right\}$.
Any solution to (\ref{eq:II}) thus satisfies $w=w_{D(w)}$ and 
therefore
\begin{eqnarray*}
D(w(t)) & = & \dps D \left( w_{D(w)}(t) \right)  \\
& = & \dps \alpha \int_{|\sigma|>1} w_{D(w)}(t,\sigma) \, d\sigma \\
& = & \dps \alpha e^{-\gamma t} \, \int_{|\sigma|>1} \left( \int_\RR 
p_0(\sigma') \varphi_{\sqrt{2
    \int_0^t D(w(s)) \, ds}}(\sigma-\sigma') \, d\sigma' \right) \,
d\sigma \\
& = & \dps e^{-\gamma t} \, F_{p_0} \left(
    \int_0^t D(w(s)) \, ds \right).
\end{eqnarray*}
It follows that the function $D(w)$ is solution in $C^0([0,+\infty[)$ 
to the
nonlinear integral equation
\begin{equation} \label{eq:inteq}
y(t)  = e^{-\gamma t}  \, F_{p_0} \left( \int_0^t y(s) \, ds \right).
\end{equation}
On the other hand, if $D \in C^0([0,+\infty[)$ is solution to
(\ref{eq:inteq}) it is easy to check that the function $w_D$ defined 
by
(\ref{eq:wD}) is solution to (\ref{eq:III}).

\medskip

\noindent
If condition (\ref{eq:condFp0}) is fulfilled, equation 
(\ref{eq:inteq})
has a unique solution in $C^0([0,+\infty[)$ (the constant function 
equal
to zero)  and the steady state $w(t,\cdot)=p_0$ thus is the unique
solution to (\ref{eq:II}) in $C^0_t(L^2_\sigma)$; otherwise, the
set of solutions to (\ref{eq:inteq}) is made of the steady state
$w(t,\cdot)=p_0$ and of the family $(y_{t_0})_{t_0 \ge 0}$ with
\begin{equation*}
y_{t_0}(t) = \, \left| \begin{array}{lll}
0 & \quad & \mbox{if } t \le  t_0  \\
z(t-t_0) & \quad & \mbox{if } t >  t_0 \end{array} \right.
\end{equation*}
where the function $z$ is defined on $[0,+\infty[$ by
\begin{equation*}
\int_0^{z(t)} \frac{dx}{F(x)} = \left\{\begin{array}{ll}
\frac{1-e^{-\gamma t}}{\gamma}\,, &\textrm{ if }\quad \gamma>0\;;\\ 
t \,,&\textrm{ otherwise\,.}
\end{array}\right.
\end{equation*}
Statement \textit{ii.}  is obtained by denoting by $v$ the solution to
(\ref{eq:III}) associated with the function $z(t)$.
\hfill$\diamondsuit$
\bigskip

\noindent\textbf{Proof of Proposition~\ref{prop:deg}:}

\medskip

\noindent
The solution $p_\epsilon$ to equation \eqref{syst:p-eps} satisfies the
inequalities
\begin{equation*}
p_\epsilon^-(t,\sigma) \le p_\epsilon(t,\sigma) \le
p_\epsilon^+(t,\sigma) \qquad \mbox{almost everywhere}
\end{equation*}
where $p_\epsilon^-$ and $p_\epsilon^+$ are defined in
$C^0_t(L^2_\sigma)$ by
\begin{equation*}
\left\{ \begin{array}{l}
\dps \partial_t p_\epsilon^- = \left(D(p_\epsilon(t))+\epsilon\right)
\partial^2_{\sigma \sigma} p_\epsilon^- -  p_\epsilon^- \, ; \\
p_\epsilon^-(0,\sigma) = p_0(\sigma) \, ; \end{array} \right. \qquad 
\qquad
\left\{ \begin{array}{l}
\dps \partial_t p_\epsilon^+ = \left(D(p_\epsilon(t))+\epsilon\right)
\partial^2_{\sigma \sigma} p_\epsilon^+ + \frac{D(p_\epsilon)}\alpha
\delta_0 \, ; \\
p_\epsilon^-(0,\sigma) = p_0(\sigma). \end{array} \right.
\end{equation*}
Therefore on the one hand
\begin{equation}
D(p_\epsilon(t)) \ge D(p_\epsilon^-(t)) = e^{-t} \, F_{p_0} \left( 
\int_0^t
  \left(D(p_\epsilon)+\epsilon\right) \right) \label{eq:binf}
\end{equation}
and on the other hand
\begin{eqnarray*}
D(p_\epsilon(t)) \le D(p_\epsilon^+(t)) & = &  F_{p_0} \left( \int_0^t
  \left(D(p_\epsilon)+\epsilon\right) \right) + \int_0^t
\frac{D(p_\epsilon)(s)}{\alpha} \, \left(\int_{|\sigma|> 1} 
\varphi_{\sqrt{2 \int_s^t
\left(D(p_\epsilon)+\epsilon\right)}} \right) \, ds  \\
& \le & F_{p_0} \left( \int_0^t
  \left(D(p_\epsilon)+\epsilon\right) \right) + \frac 1 \alpha 
\int_0^t
D(p_\epsilon)(s) \, ds.
\end{eqnarray*}
If (\ref{eq:condFp0}) is
not fulfilled, using (\ref{eq:binf}) and the property that $F_{p_0}$
is strictly increasing on $[0,+\infty[$, we obtain that
\begin{equation*}
D(p_\epsilon) \ge z(t)
\end{equation*}
where $z(t)$ is the function defined in the proof of 
Lemma~\ref{lem:2}. As
for any $0< t_0 \le T$, there exists $\eta > 0$ such that $z(t) \ge
\eta$ on $[t_0,T]$ the same reasoning as in the non-degenerate case 
leads to the conclusion
that $(p_\epsilon)$ converges up to an extraction to $p$ in
${\cal D}'(]0,+\infty[\times \RR)$ and in
$L^2([t_0,T],L^2(\RR))$ for any $0 < t_0 < T < +\infty$, $p$
being a solution to \eqref{eq:syst-p} in $C^0(]0,+\infty[,L^2_\sigma)$
such that $D(p) > 0$ on $]0,+\infty[$. \hfill$\diamondsuit$
\medskip

\vskip10pt
\section{Steady states}\label{sec:stat}
Throughout this section the shear rate $b$ is assumed to be a given 
constant and we are looking
for solutions in $L^1(\RR)$ to the following system~:
\begin{EqSystem}\label{eq:HLstat}
 -b\; \partial_{\sigma}p+D(p)\;
    \partial^2_{\sigma\sigma}p-\chi_{\RR\setminus [-1,1]}
    \;p+\frac{D(p)}{\alpha}\;\delta_{0}(\sigma)=0\quad \textrm{ on } 
(0;T)\times \RR\;;\label{eq:pstat}
\\
p\,\ge\,0\;,\int_\RR p=1\;;\\
D(p) = \alpha\, \int_{|\sigma|>1} p(\sigma)\, d\sigma\;.
\end{EqSystem}
Our main results are summarized in 
the following~:
\begin{proposition}\label{prop:solstat}{\tiny .}
\begin{itemize}
\item[i] If $b\equiv 0$, any probability density which is compactly
  supported in $[-1;+1]$ is a solution  
to~\eqref{eq:pstat} which satisfies $D(p)=0$. If  
$\alpha\leq\frac{1}2$,
these are the only stationary  
solutions (and there are infinitely many), whereas when
$\alpha>\frac{1}2$ there exists a unique stationary solution
corresponding to  
 a positive value of  $D$, which is explicitly  given
 by~\eqref{eq:solstat0} and \eqref{eq:D0} 
below. This solution is even and with exponential decay at infinity.
 \item[ii] If $b\not\equiv 0$, for any $\alpha > 0$,  there exists  a 
 unique  stationary solution 
   to~\eqref{eq:pstat}, and it corresponds to  
 a positive value for  $D$, which is implicitly given
 by~\eqref{eq:solstat} and \eqref{eq:norm-new} below. This solution 
has
 exponential decay at infinity.
\end{itemize}
\end{proposition}
\begin{remarque} The statement in the above proposition is already 
pointed 
out by H\'ebraud and
Lequeux~\cite{bi:HL}.
\end{remarque}

\vskip6pt\noindent\textbf{Proof of Proposition~\ref{prop:solstat}:}\\

\noindent\textbf{The case when $b\equiv 0$}\\
We first observe that any non-negative function $p$ which is 
normalized in  $L^1(\RR)$ and with
support in $[-1;+1]$ is a solution to the system~\eqref{eq:HLstat} 
since in that case
 all terms in equation~\eqref{eq:pstat} cancel.  We now
 examine the issue of
existence of solutions of~\eqref{eq:HLstat} such that $D(p)>0$. For 
simplicity we denote
$D=D(p)$. For given constant $D>0$, it is  very easy to calculate 
explicitly the solutions
of~\eqref{eq:pstat} on each of the three regions $\sigma<-1$,
$\sigma\in\,[-1;+1]$ 
and $\sigma>1$. Using compatibility conditions on $\RR$ and the fact 
that $p$ has to be in
$L^1(\RR)$ one  obtains~: 
\begin{equation}\label{eq:solstat0}
p(\sigma)=\left\{
\begin{array}{rl}
\dps \frac{\sqrt{D}}{2\,\alpha}\,e^{\,(1+\sigma)/\sqrt{D}}&\quad
\textrm{ if }\sigma\leq -1\;, \vspace{3mm} \\ \vspace{3mm}
\dps \frac{1}{2\,\alpha}\,\sigma+\frac{\sqrt{D}+1}{2\,\alpha}&\quad 
\textrm{ if }-1\leq \sigma\leq 0\;,\\ \vspace{3mm}
\dps -\frac{1}{2\,\alpha}\,\sigma+\frac{\sqrt{D}+1}{2\,\alpha}&\quad 
\textrm{ if }0\leq \sigma\leq 1\;,\\ 
\dps \frac{\sqrt{D}}{2\,\alpha}\,e^{(1-\,\sigma)/\sqrt{D}}&\quad 
\textrm{ if }1\leq \sigma\;.
\end{array}
\right.
\end{equation}  
The compatibility condition $D=D(p)$ happens to be then automatically 
satisfied and the normalization 
constraint     $\int_{\RR}p=1$ 
imposes that $D$ 
solves 
\begin{equation}\label{eq:D}
D+\sqrt{D}=\alpha-\frac{1}2\;.
\end{equation}
Since $D\geq 0$, we immediately reach a contradiction when 
$\alpha<\frac{1}2$, whereas when $\alpha> \frac{1}2$ 
equation~\eqref{eq:D} 
admits a unique positive solution; namely 
\begin{equation}\label{eq:D0}
D=-\frac{1}2+\frac{\sqrt{4\,\alpha}-1}2\;.
\end{equation}

\noindent\textbf{The case when $b\not\equiv 0$}\\
First of all, we observe that if $D=0$ every term in
equation~\eqref{eq:pstat} but $b\,\partial_\sigma p$  vanish. Thus $p$
has to be a non-zero constant which is in contradiction with $p\in 
L^1(\RR)$. So necessarily $D>0$. For given positive constant $D$, we 
then  
solve~\eqref{eq:pstat} as above and obtain
\begin{equation}\label{eq:solstat}
p(\sigma)=\left\{
\begin{array}{rl}
 \dps a_1\,e^{\,\beta^+\,\sigma}&\quad \textrm{ if }\sigma\leq -1\;,
 \vspace{3mm} \\
\dps a_2 e^{\frac{b}{D}\,\sigma} + a_2 - \frac D {b\alpha} &\quad 
\textrm{ if
}-1\leq \sigma\leq 0\;, \vspace{3mm} \\ 
\dps \big(a_2 - \frac D {b\alpha}\big) e^{\frac{b}{D}\,\sigma}  + a_2 
&\quad
\textrm{ if }0\leq
\sigma\leq 1\;, \vspace{3mm} \\
\dps a_1\,e^{\,\beta^-\,\sigma}&\quad \textrm{ if }1\leq \sigma\;,
\end{array}
\right.
\end{equation}  
with 
\begin{equation*}
\beta^\pm=\frac{b}{2\,D}\pm\frac{1}{2}\,\sqrt{\frac{b^2+4\,D}{D^2}}\,,
\end{equation*}
\begin{equation*}
a_1=\frac{e^{\frac 1 2  \sqrt{\frac{b^2}{D^2}+\frac{4}{D}}}}
{\alpha \left( \beta^+ e^{b/2D} - \beta^- e^{-b/2D} \right)},  
\end{equation*}
and
\begin{equation*}
a_2=\frac{D \beta^+ e^{b/2D}}{\alpha b \left( \beta^+ e^{b/2D} - 
\beta^-
    e^{-b/2D} \right)}. 
\end{equation*}
It is tedious but easy to check that this function always fulfills the
self-consistency condition $D=D(p)$ 
and that the normalization condition $\int_\RR p=1$ reads 
\begin{equation}\label{eq:norm-new}
\frac{D}{b}\,\frac{ (1 + \beta^+)
    +( \beta^- -1)\,e^{-b/D}} {\beta^+-\beta^-\,e^{-b/D}} +D=\alpha\;.
\end{equation} 
{}For any $b>0$ (the negative values of
$b$ are dealt with by replacing $\sigma$ by $-\sigma$), the left-hand
side of (\ref{eq:norm-new}) is a continuous function which goes to 
$+\infty$ when $D$ goes to infinity  and goes to zero  when $D$ goes
to $0$. This already ensures the existence of at least one steady 
state for any $\alpha > 0$. Moreover,  setting $z=\frac{b^{2}}{D}$ (for example) 
we may rewrite the left-hand side of~(\ref{eq:norm-new}) as 
\ben
f(z) = \frac{b^2}{z} + \frac{2\,b^{2}}{z}\, 
\left[ 
\frac{ 1+ \frac 1 {2b} z \coth(z/2b) +\frac 1 {2b} (z^2+4z)^{1/2}}
{ z + (z^2+4z)^{1/2} \coth(z/2b) }
\right].
\een
Next we check that the function $f$ is monotone decreasing (thus, 
the left-hand side of~\eqref{eq:norm-new} is increasing with respect to $D$), 
whence the uniqueness result. 
\hfill$\diamondsuit$

\vskip10pt

\noindent\textbf{Acknowledgements.} We would like
to thank Philippe Coussot for pointing out the H\'ebraud-Lequeux 
equation
to us. We also warmly  thank Claude Le Bris for stimulating 
discussions.

\bibliographystyle{plain}

\end{document}